\documentclass[10pt]{amsart}
\usepackage{amsmath}
\usepackage[curve,matrix,arrow]{xy}
\title[Contact Cohomology of the Projective Plane]
{Contact Cohomology of the Projective Plane}
\author{Lars Ernstr{\"o}m }
\address{ Department of Mathematics,
Royal Institute of Technology,
S-100 44 Stockholm, Sweden} 
\email{ernstrom@math.kth.se}
\author{Gary Kennedy}
\address{Ohio State University at Mansfield, 1680 Univ. Dr.,
Mansfield, Ohio 44906, USA}
\email{kennedy@math.ohio-state.edu}
\thanks{The first author was supported by the Swedish National Research Council.}
\keywords{Characteristic number, Deligne-Mumford stable curve,
Kontsevich-stable map, Gromov-Witten invariant, quantum cohomology.}
\subjclass{Primary 14C17, 14N10. Secondary 14D22.}

\newcommand{\abstracttext}
{We construct an associative ring which is a deformation of the quantum cohomology ring of the projective plane. Just as the quantum cohomology encodes the incidence characteristic numbers of rational plane curves, the contact cohomology encodes the tangency characteristic numbers.}
\numberwithin{equation}{section}
\newtheorem{thm}[equation]{Theorem}
\newtheorem{prop}[equation]{Proposition} 

\newtheorem{lem}[equation]{Lemma}

\theoremstyle{definition}

\theoremstyle{remark}

\newcommand{\q}{\mathbf Q} 
\newcommand{\pl}{{\mathbf P^1}} 
\newcommand{\pp}{{\mathbf P^2}} 
\newcommand{\dpp}{\mathbf {\check P}^{2}} 
\newcommand{\ap}{A^*(\pp)} 
\newcommand{\ai}{A^*(I)} 
\newcommand{\qt}{\q[[T]]}
\newcommand{\apt}{\ap[[T]]} 
\newcommand{\ait}{\ai[[T]]}%
\newcommand{\hd}{\check{h}} 
\newcommand{\dd}{\check{d}} %
\newcommand{\stablefour}{{\overline M}_{0,4}} 
\newcommand{\pmoduli}[2]{{\overline M}_{0,#1}(\pp,#2)}
\newcommand{\lamoduli}[2]{{\overline M}^1_{0,#1}(\pp,#2)}
\newcommand{\imoduli}[2]{{\overline M}_{0,#1}(I,#2)}
\newcommand{\double}[1]{{\overline M}_{0,#1,2}(I,(0,2))}
\newcommand{\spmoduli}{{\overline M}(m,n,d)}
\newcommand{\spamoduli}[3]{{\overline M}(#1,#2,#3)}
\newcommand{\Mone}{M_1}
\newcommand{\Mtwo}{M_2}
\newcommand{\kp}[3]
{\frac{\partial^3 \calK}{\partial y_{#1}\partial y_{#2}\partial z_{#3}}}

\newcommand{\calN}{{\mathcal N}}

\newcommand{\calR}{{\mathcal R}}

\newcommand{\calK}{{\mathcal K}}

\begin{document}
\begin{abstract} \abstracttext \end{abstract} \maketitle
%
%
%
%
\section{Introduction}
\label{intro}
\par
In this paper we construct an associative ring which we call the {\em contact cohomology ring} of the projective plane. We believe that an analogous construction should work for all homogeneous varieties, but in our proof of associativity we rely on certain technical results from our earlier paper on characteristic numbers of rational plane curves \cite{ErnstromKennedy}.
\par
As we formulate it in section \ref{qtp} of the present paper, the quantum product is actually a whole family of products parametrized by elements of the Chow ring $\ap$. Each product is defined on the formal power series ring $\apt$ in one variable, and encodes the characteristic numbers
\begin{equation*}
N_d=
\text{ the number of rational plane curves of degree $d$
through $3d-1$ general points}.
\end{equation*}
The remarkable fact about the quantum product is that it is associative. This fact, together with the triviality $N_1=1$, suffices to determine all values of $N_d$.
\par
The contact products are defined on the same formal power series ring. But now the parameter space is the Chow ring $A^*I$ of the incidence variety of points and lines in $\pp$, and these products 
encode (as we explain in section \ref{rrc}) a larger collection of characteristic numbers:
\begin{align*}
N_d(a,b,c)=&\text{ the number of rational plane curves of degree $d$
through $a$ general points, tangent to}\\ &\text{$b$ general lines, and
tangent to $c$ general lines at a specified general point on each line}\\
&\quad\text{(where $a+b+2c=3d-1$).}
\end{align*}
Section \ref{qcp} is devoted to defining these products precisely and to showing that they are likewise associative. As we explain in section \ref{rrc}, the associativity implies a remarkable recursive relation among the characteristic numbers, but it is insufficient to determine all their values unless one already knows all values $N_d(a,b,c)$ for which $a<3$. (Our previous paper \cite{ErnstromKennedy} explains how to obtain this additional information.)
\par
The associativity of the quantum product is a consequence of the recursive structure of the boundary divisors on the moduli space of stable maps to $\pp$ (of a given degree and with a given number of markings). In brief, each boundary divisor is a fiber product of simpler instances of the same sort of moduli space. But in studying questions of tangency it is natural to look at a moduli space of stable lifts (defined in section \ref{siv}), whose boundary divisors have a somewhat more complicated structure. To understand these divisors, consider a family of immersions $\pl \to \pp$. Associated to each immersion is its lift, a map from $\pl$ to the incidence correspondence $I$. Now suppose that the family of immersions degenerates to a map from a two-component curve onto the union of two plane curves of lower degree. Then one can show that the family of lifts degenerates to a map from a three-component curve, with the central component mapping two-to-one onto the fiber of $I$ over a point of $\pp$. Thus to create a moduli space with a recursive structure we form a fiber product of the space of stable lifts with a space of two-to-one covers of fibers of $I$ over $\pp$; we then exploit this recursive structure to define our associative product.
\par
Our recursive relation specializes to that of Di~Francesco and Itzykson
\cite{FrancescoItzykson}, who also interpret their formula as the
associativity of a certain product. We thank S.~Colley for carefully reading earlier versions of this paper and making many helpful suggestions.
We also want to indicate our
indebtedness to B.~Fantechi, W.~Fulton, L.~G{\"o}ttsche and R.~Pandharipande.
\par
We should warn the reader of a potentially confusing clash in terminology: in symplectic geometry Y.~Eliashberg and H.~Hofer have introduced a new invariant of contact manifolds called ``contact homology.'' The Gromov-Witten invariants of a symplectic manifold with contact boundary take values in the contact homology of the boundary.
\par
\section{A ``push-pull'' formula}
\label{ppf}
\par
For an algebraic variety, scheme, or algebraic stack $X$, we write $A_*X$ for the rational equivalence group with coefficients in $\q$, and call it the {\em homology}; we write $A^*X$ for the operational cohomology ring. (See \cite{Fulton} for the case of a variety or scheme, \cite{Vistoli} for intersection theory on stacks.) We will often use the fact that $A_*(X)$ and
$A^*(X)$ of a nonsingular variety $X$ are naturally isomorphic. 

We record here a formula which we will use repeatedly. It is the stack version of \cite[Proposition 1.7]{Fulton}.
\begin{lem}
\label{pushpull}
Suppose that
$$ \xymatrix{
& M_1\times_X M_2 \ar[dl]^{q_1} \ar[dr]_{q_2} \\
M_1 \ar[dr]^{g_1} & & M_2 \ar[dl]_{g_2} \\
& X }
$$
is a fiber square of stacks over a nonsingular variety $X$, with $g_2$ flat and $g_1$ proper. Then $q_1$ is flat, $q_2$ is proper, and for each class $\alpha \in A^*M_1$ we have
$$
q_{2*}\left(q_1^*\alpha\cap\lbrack M_1 \times_X M_2\rbrack\right)
=g_2^*g_{1*}(\alpha\cap\lbrack M_1 \rbrack)\cap\lbrack M_2 \rbrack.
$$
\end{lem}
\par
\section{The quantum product}
\label{qtp}
Here we recall the basic definitions of quantum cohomology and prove that the quantum product on $\pp$ is associative. Our proof is basically the same as in section 8 of \cite{FultonP}, except that we avoid the use of coordinates; we use this proof as a prototype for the longer proof in section \ref{qcp}. As general references for the material in this section we suggest \cite{BehrendManin}, \cite{ErnstromKennedy}, \cite{FultonP}, \cite{Kontsevich}, \cite{KontsevichManin}, and \cite{LiTian}.
\par
Let $\pmoduli{n+3}{d}$ denote the stack of stable maps from curves of arithmetic genus $0$, with $n+3$ markings, to the projective plane. Let $i,j,k:\pmoduli{n+3}{d}\to \pp$ be the evaluation maps associated to the first three markings; let $e_1,\dots,e_n$ be the others. Let $\ap$ be the Chow cohomology ring of the plane with rational coefficients, and let $\apt$ be the ring of formal power series in one variable. Suppose that $\alpha$, $\beta$, and $\delta$ are elements of $\ap$. Then the {\em quantum product of $\alpha$ and $\beta$, deformed by $\delta$}, is the element of $\apt$ whose $n$th coefficient is
\begin{equation}
\label{qtpn}
\frac{1}{n!}
\sum_{d \geq 0}
 k_*(i^*\alpha \cup j^*\beta \cup
{\bigcup_{t=1}^n} e_t^*\delta \cap
\lbrack \pmoduli{n+3}{d} \rbrack ).
\end{equation}
Note that, since the dimension of $\pmoduli{n+3}{d}$ is $3d+2+n$, the sum is finite. We will denote the quantum product by $(\alpha * \beta)_{\delta}$ or simply $\alpha * \beta$. Extending by $\qt$-linearity, we have a product on $\apt$. We call this ring the {\em quantum cohomology} of $\pp$, and denote it
by $QH^*(\pp)$.
\begin{prop}
For each $\delta$ in $\ap$, the quantum product is commutative and associative. The identity element $1 \in \ap$ for the ordinary cup product also serves as the identity element for the quantum product.
\end{prop}
\begin{proof}
The commutativity is obvious. Unless $n=d=0$ the class
$$
i^*1 \cup j^*\beta \cup
{\bigcup_{t=1}^n} e_t^*\delta
$$
is the pullback via the forgetful morphism $\pmoduli{n+3}{d} \to \pmoduli{n+2}{d}$ of a class on the latter space. Since the fibers of this morphism have positive dimension, the projection formula tells us that the corresponding term of (\ref{qtpn}) vanishes. As for the remaining term $n=d=0$, note that in this case $i$, $j$, and $k$ are all the same isomorphism from $\pmoduli{3}{0}$ to $\pp$. Hence $1$ is the identity element for the quantum product.
\par
To see that the quantum product is associative, consider the moduli space $\pmoduli{n+4}{d}$. Denote the first four evaluation maps $\pmoduli{n+4}{d} \to \pp$ by $i$, $j$, $k$, and $l$, and the remaining $n$ maps by $e_t$. Consider the ``forgetful'' morphism from $\pmoduli{n+4}{d}$ to $\stablefour$, the space of stable 4-pointed curves, which associates to a stable map its source curve together with the first four markings, with any unstable components contracted to a point.
The space $\stablefour$ is isomorphic to $\pl$. It has a distinguished point $P(12\mid 34)$ representing the two-component curve having the first two markings on one component and the latter two on the other; similarly there are two other distinguished points $P(13\mid 24)$ and $P(14\mid 23)$. Hence on $\pmoduli{n+4}{d}$ there are three linearly equivalent divisors $D(12\mid 34)$, $D(13\mid 24)$, and $D(14\mid 23)$, which we call the {\em special boundary divisors}.
\par
Kontsevich \cite{Kontsevich} (cf. \cite{FultonP}) identifies the components of $D(12\mid 34)$. For a finite set $A$, let $\pmoduli{A}{d}$ denote the stack of stable maps with markings labeled by $A$. Suppose that $A_1 \cup A_2$ is a partition of $\{1,\dots,n+4\}$ and that $d_1+d_2=d$. Suppose that $\{\star\}$ is a single-element set. Then the fiber product
$$
D(A_1,A_2 ,d_1,d_2) =
\pmoduli{A_1\cup\{\star\}}{d_1} \times_{\pp}
\pmoduli{A_2\cup\{\star\}}{d_2}
$$
is naturally a substack of $\pmoduli{n+4}{d}$; the typical point represents a map from a curve with two components, with the point of attachment corresponding to the point labeled by $\star$, as indicated in Figure 3.1.
The divisor $D(12\mid 34)$ is the sum
\begin{equation}
\label{donetwo}
D(12\mid 34)=\sum D(A_1,A_2 ,d_1,d_2)
\end{equation}
over all partitions $A_1 \cup A_2$ in which $1$ and $2$ belong to $A_1$, and $3$ and $4$ belong to $A_2$. There are corresponding statements for the other two special boundary divisors.
\par
%
%
$$
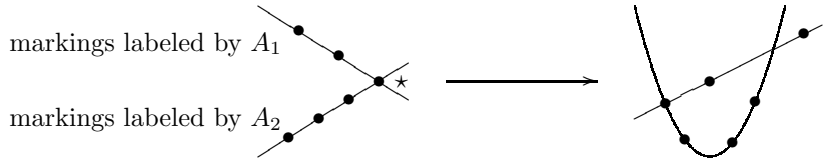
\xy
0;<1cm,0cm>:
(-1,-1);(1,0.25)**\dir{-};
(-1,1);(1,-0.25)**\dir{-};
(0.9,0)*{\star};
(0.6,0)*{\bullet};
(-0.467,0.667)*{\bullet};
(0.067,0.333)*{\bullet};
(-0.6,-0.75)*{\bullet};
(-0.2,-0.5)*{\bullet};
(0.2,-0.25)*{\bullet};
(-2.5,0.5)*{\text{markings labeled by }A_1};
(-2.5,-0.5)*{\text{markings labeled by }A_2};
(1.5,0);(3.5,0)**\dir{-};
?>*\dir{>};
(4,-0.5);(6.5,0.75)**\dir{-};
(5,0)*{\bullet};
(6.25,0.625)*{\bullet};
(4,1);(6,1)**\crv{(5,-3)}\POS?(0.333)*{\bullet}\POS?(0.65)*{\bullet}\POS?(0.8)*{\bullet};
(4.41,-0.3)*{\bullet};
\endxy
$$
\begin{center}
{\bf Figure 3.1} Typical member of $D(A_1,A_2 ,d_1,d_2)$.
\end{center}
\smallskip
\par
Since the divisors
$D(12\mid 34)$ and $D(14\mid 23)$ are linearly equivalent, the equation
\begin{multline}
\label{transeq}
\frac{1}{n!}
\sum_{d \geq 0}
l_*(i^*\alpha \cup j^*\beta \cup k^*\gamma \cup {\bigcup_{t=1}^n}
e_t^*\delta \cap \lbrack D(12\mid 34) \rbrack ) \\
= \frac{1}{n!}
\sum_{d \geq 0}
l_*(i^*\alpha \cup j^*\beta \cup k^*\gamma \cup {\bigcup_{t=1}^n}
e_t^*\delta \cap \lbrack D(14\mid 23) \rbrack ).
\end{multline}
is valid for each triple $\alpha$, $\beta$, $\gamma$ of classes in $\ap$ and for each $n$. According to (\ref{donetwo}), the divisor $D(12\mid 34)$ is a sum of fiber products, each of which fits into a fiber diagram
$$ \xymatrix{
& & \pmoduli{n_1+3}{d_1}\times_{\pp}\pmoduli{n_2+3}{d_2}
\ar[ddl]^{q_1} \ar[ddr]_{q_2} & \\
\\
& \pmoduli{n_1+3}{d_1} \ar@/^/[dd]^{e_{1t}}
\ar@/_/[dd]_{i_1}\ar[dd]|{j_1} \ar[ddr]^{g_1} & &
\pmoduli{n_2+3}{d_2} \ar@/^/[dd]^{e_{2t}} \ar@/_/[dd]_{k_2}
\ar[dd]|{l_2} \ar[ddl]_{g_2} & \\
\\
& \pp & \pp & \pp & }
$$
in which $i=i_1 \circ q_1$, $j=j_1 \circ q_1$, $k=k_2 \circ q_2$, and $l=l_2 \circ q_2$. Furthermore each $e_t$ equals either $e_{1t} \circ q_1$ or $e_{2t} \circ q_2$; by relabeling we may assume that the former equation holds for $1 \leq t \leq n_1$. Note that, for specified partitions $n=n_1+n_2$ and $d=d_1+d_2$, the number of such fiber diagrams is $\binom{n}{n_1}$. Let $\Mone=\pmoduli{n_1+3}{d_1}$ and $\Mtwo=\pmoduli{n_2+3}{d_2}$. By Lemma~\ref{pushpull} and the projection formula we have
\begin{align*}
l_*&\left(i^*\alpha\cup j^*\beta\cup k^*\gamma\cup\bigcup_{t=1}^n
e_t^*\delta\cap[\Mone\times_{\pp}\Mtwo]\right) \\
 &=l_{2*}q_{2*}\left(q_1^*(i_1^*\alpha\cup
j_1^*\beta\cup\bigcup_{t=1}^{n_1}e_{1t}^*\delta)\cup
q_2^*(k_2^*\gamma\cup\bigcup_{t=n_1+1}^{n}e_{2t}^*\delta)
\cap[\Mone\times_ {\pp}\Mtwo]\right) \\
&=l_{2*}\left(
\left(k_2^*\gamma\cup\bigcup_{t=n_1+1}^n e_{2t}^*\delta\right)
\cap q_{2*}
\left(q_1^*(i_1^*\alpha\cup
j_1^*\beta\cup\bigcup_{t=1}^{n_1}e_{1t}^*\delta)
\cap[\Mone\times_{\pp}\Mtwo] \right)
\right) \\
&=l_{2*}\left(
\left(k_2^*\gamma\cup\bigcup_{t=n_1+1}^ne_{2t}^*\delta\right)
\cap
\left(g_2^*g_{1*}(i_1^*\alpha\cup
j_1^*\beta\cup\bigcup_{t=1}^{n_1}e_{1t}^*\delta
\cap[\Mone])\cap[\Mtwo]
\right)
\right) \\
&=l_{2*}\left(
g_2^*\left( g_{1*}(i_1^*\alpha\cup
j_1^*\beta\cup\bigcup_{t=1}^{n_1}e_{1t}^*\delta\cap[\Mone])
\right) \cup
k_2^*\gamma\cup\bigcup_{t=n_1+1}^ne_{2t}^*\delta
\cap[\Mtwo]
\right).
\end{align*}
Summing---for fixed $n$---over all $d \geq 0$ and all components of $D(12\mid 34)$, we find that the left side of (\ref{transeq}) equals the $n$th coefficient of $(\alpha*\beta) *\gamma$. A similar argument shows that the right side of (\ref{transeq}) equals the $n$th coefficient of $\alpha*(\beta*\gamma)$.
\par
Invoking the $\qt$-linearity, we conclude that the quantum product is associative.
\end{proof}
\section{Stable maps to the incidence variety}
\label{siv}
\par
Let $I$ be the incidence variety of points and lines in $\pp$. Its cohomology ring $A^*(I)$ is generated by two classes $h$ and $\hd$ representing the pullbacks of, respectively, the class of a line in $\pp$ and the class of a dual line in $\dpp$. The fundamental class of a curve is determined by its intersection numbers $d$ and $\dd$ with the two classes; we denote this class by $(d,\dd)$. The incidence variety is a projective line bundle over the plane; we denote the projection map by $p:I \to \pp$.
For sake of simplicity in notation, we will consider $\ap$ as a subring
of $A^*(I)$, embedded by
$$
p^*\colon \ap \hookrightarrow A^*(I).
$$
\begin{lem} \label{lemmachar}
$
\ap=\{\alpha \in A^*(I) \mid p_*\alpha=0\}.
$
\end{lem}
(Note that whenever we use a push-forward of a cohomology class on a nonsingular
space, it is defined via the push-forward of the dual homology class.)
\par
As we explain in our earlier paper \cite{ErnstromKennedy}, an $n$-pointed immersion from $\pl$ to $\pp$ can be lifted to an $n$-pointed map from $\pl$ to $I$; if the immersion has degree $d$ then the class of the lift is $(d,2d-2))$. Thus there is a substack $M^1_{0,n}(\pp,d)$ of $\imoduli{n}{(d,2d-2))}$ representing the lifts of immersions. We call its closure $\lamoduli{n}{d}$ the {\em stack of stable lifts}; it is birationally isomorphic to the stack $\pmoduli{n}{d}$ and thus has dimension $3d-1+n$.
\par
We will be working with two other special stacks of stable maps to the incidence variety. For $n\geq2$, consider the stack ${\overline M}_{0,n}(I,(0,1))$, whose typical member is an $n$-pointed isomorphism from $\pl$ to a fiber of $I$ over $\pp$. Let $e_1,\dots,e_n$ be the evaluation maps.
\begin{lem} \label{disappear}
Let $\alpha$ be a class in $\ap$ and let $\delta_2,\dots,\delta_n$ be classes in $A^*(I)$. Then
$$
e_{ 2*}\left(e_1^*p^*\alpha\cup e_2^*\delta_2\cup \dots 
\cup e_n^*\delta_n\cap [{\overline M}_{0,n}(I,(0,1))]\right) =0.
$$
\end{lem}
\begin{proof}
Since each element of ${\overline M}_{0,n}(I,(0,1))$ is a map to a fiber of $I$ over $\pp$, we have that $p\circ e_1=p\circ e_2$. Thus $e_1^*p^*=e_2^*p^*$, and by the projection formula
\begin{align*}
&e_{ 2*}\left(e_2^*p^*\alpha\cup e_2^*\delta_2\cup \dots 
\cup e_n^*\delta_n\cap [{\overline M}_{0,n}(I,(0,1))] \right) =\\
&p^*\alpha\cap e_{2*}\bigl(e_2^*\delta_2\cup\dots\cup e_n^*\delta_n
\cap[{\overline M}_{0,n}(I,(0,1))] \bigr).
\end{align*}
The class in parentheses is the pullback of a class on ${\overline M}_{0,n-1}(I,(0,1))$
via the morphism
$$
{\overline M}_{0,n}(I,(0,1))\to {\overline M}_{0,n-1}(I,(0,1))
$$
which forgets the first point. Furthermore $e_2$ factors through this forgetful morphism. Hence by the projection formula its pushforward via $e_2$ vanishes.
\end{proof}
\par
Similarly, let $M_{0,n,2}(I,(0,2))$ be the substack of $\imoduli{n+2}{(0,2)}$ representing degree $2$ maps from $\pl$ to a fiber of $I$ over $\pp$, with the two ramification points specially marked and with an additional $n$ markings. Let $\double{n}$ be its closure, a stack of dimension $4+n$.
Denote by
$e_1,\dots,e_{n+2}$ the evaluation maps of $\double{n}$. It does not matter
which two are evaluation maps at the ramification; this flexibility is
convenient in proving the next three lemmas. (Later we will use a different
notation.)
\begin{lem} \label{lemmafundclass}
Let $\delta_2,\dots,\delta_{n+2}$ be classes in $A^*(I)$.
Then
$$
e_{2*}\bigl(e_2^*\delta_2\cup\dots\cup e_{n+2}^*\delta_{n+2}
\cap[\double{n}]\bigr)=0.
$$
\end{lem}
(Note that there is no condition on the marking corresponding to $e_1$.)
\begin{proof}
Each evaluation map $e_i$ ($i=2,\dots,n+2$) can be factored into the product map $f=e_2\times\dots\times e_{n+2}$ followed by projection $f_i$ onto the appropriate factor.
$$ \xymatrix{
\double{n} \ar[d]^{f} \\
I\times_\pp\dots \times_\pp I
\ar[d]^{f_i} \\
I}
$$
Therefore
\begin{align*}
&e_{2*}\bigl(e_2^*\delta_2\cup\dots\cup e_{n+2}^*\delta_{n+2}
\cap[\double{n}]\bigr)=\\
&f_{2*}f_*\bigl(f^*(f_2^*\delta_2\cup\dots\cup f_{n+2}^*\delta_{n+2})
\cap[\double{n}]\bigr)=\\
&f_{2*}\bigl(f_2^*\delta_2\cup\dots\cup f_{n+2}^*\delta_{n+2}
\cap f_*[\double{n}]\bigr),
\end{align*}
which is zero because $\dim(\double{n})=4+n$ and
$\dim(I\times_\pp\dots \times_\pp I)=3+n$.
\end{proof}
\begin{lem}\label{lemmapbpp}
Let $\alpha$ be a class in $\ap$ and let $\delta_2,\dots,\delta_{n+2}$ be classes in $A^*(I)$. Then
$$
e_{2*}\left(e_1^*p^*\alpha\cup e_2^*\delta_2\cup \dots 
\cup e_{n+2}^*\delta_{n+2}\cap [\double{n}]\right)=0.
$$
\end{lem}
\begin{proof}
Since each element of $\double{n}$ is a map to a fiber of $I$ over $\pp$, we have that $p\circ e_1=p\circ e_2$. Thus $e_1^*p^*=e_2^*p^*$, and by the projection formula
\begin{align*}
&e_{2*}\left(e_{2}^*p^*\alpha\cup e_2^*\delta_2\cup \dots 
\cup e_{n+2}^*\delta_{n+2}\cap [\double{n}]\right)=\\
&p^*\alpha\cap e_{2*}\bigl(e_2^*\delta_2\cup\dots\cup e_{n+2}^*\delta_{n+2}
\cap [\double{n}]\bigr),
\end{align*}
which is zero by Lemma~\ref{lemmafundclass}.
\end{proof}
\begin{lem}\label{lemmapfp}
For any classes $\delta_2,\dots,\delta_{n+2}$ in $A^*(I)$, the class
$$
e_{1*}\left(e_2^*\delta_2\cup \dots \cup e_{n+2}^*\delta_{n+2}\cap [\double{n}]\right)
$$
is in $\ap$.
\end{lem}
\begin{proof}
Since $p\circ e_1=p\circ e_2$, we have
\begin{align*}
&p_*e_{1*}\left(e_2^*\delta_2\cup \dots \cup e_{n+2}^*\delta_{n+2}\cap [\double{n}]\right)=\\
&p_*e_{2*}\bigl(e_2^*\delta_2\cup \dots \cup e_{n+2}^*\delta_{n+2}\cap [\double{n}],
\end{align*}
which is zero by Lemma~\ref{lemmafundclass}. Now apply Lemma~\ref{lemmachar}.
\end{proof}
\par
\section{The contact product}
\label{qcp}
\par
Denote by $\spmoduli$ the fiber product
$$
\lamoduli{m+3}{d} \times_I \double{n}.
$$
Denote the evaluation maps $\spmoduli \to I$ coming from $\lamoduli{m+3}{d}$ by $i,j,e_1,\dots,e_{m}$ and $e_{\star}$; denote those 
coming from $\double{n}$ by $f_1,\dots,f_{n},f_{\star}$ and $k$, with the latter two being at the points of ramification. The fiber product over $I$ is defined using the maps $e_{\star}$ and $f_{\star}$. Thus the typical point of $\spmoduli$ represents a map from a curve with two components, as depicted in Figure 5.1; the vertical component maps two-to-one to a fiber of $I$ and is ramified at the points marked $\times$. In addition to the indicated markings, there are $m$ additional markings on the horizontal component and $n$ on the vertical component.
%
%
%
%
$$
\xy
0;<1cm,0cm>:
(0,1);(2.5,1)**\dir{-};
(2,0);(2,1.25)**\dir{-};
(0.5,1)*{\bullet}+U(4)*{i};
(1,1)*{\bullet}+U(4)*{j};
(2,1)*{\bullet}+UL(2)*{e_\star};
(2,1)*{\bullet}*{\times}+DR(2.5)*{f_\star};
(2,0.25)*{\bullet}*{\times}+R(2)*{k};
\endxy
$$
\begin{center}
{\bf Figure 5.1} Typical member of $\spmoduli$.
\end{center}
\smallskip
\par

Given cohomology classes $\alpha$ and $\beta$ in $\ap$ and a class $\delta$ in $A^*(I)$, let $(\alpha \bullet \beta)_\delta$ (or simply $\alpha \bullet \beta$) be the element of $\ait$ whose $q$th coefficient is
\begin{equation}\label{contactproduct}
\sum_{\substack{m + n = q\\d>0}}
\frac{2}{m!n!}
k_*
\bigl(i^*\alpha\cup j^*\beta\cup\bigcup_{t=1}^m
e_t^*\delta\cup\bigcup_{t=1}^nf_t^*\delta\cap [\spmoduli]\bigr).
\end{equation}
Note that with our convention $\ap$ is a subring of $\ai$, so it makes sense to pull back classes of $\ap$ via evaluation maps to $I$. The product formula (\ref{contactproduct}) makes sense for any two classes $\alpha$ and $\beta$ in $A^*(I)$. However, to prove that the product is associative we will need to assume that the classes are in $\ap$.
\par
We now prove that $\alpha \bullet \beta$ is in $\apt$. (For this proof we don't need to assume that $\alpha$ and $\beta$ are in $\ap$.) Let $i_1$, $j_1$, $e_{1\star}$ and $e_{1t}$ indicate the evaluation maps from $\lamoduli{m+3}{d}$ to $I$; thus if $q_1$ is the projection of $\spmoduli$ onto its first factor, we have $i=i_1 \circ q_1$, etc. Similarly let $q_2$ be the projection onto the second factor and 
$k_2,f_{2\star},f_{2t}:\double{n} \to I$ be the evaluation maps, so that $k=k_2 \circ q_2$, etc.
$$ 
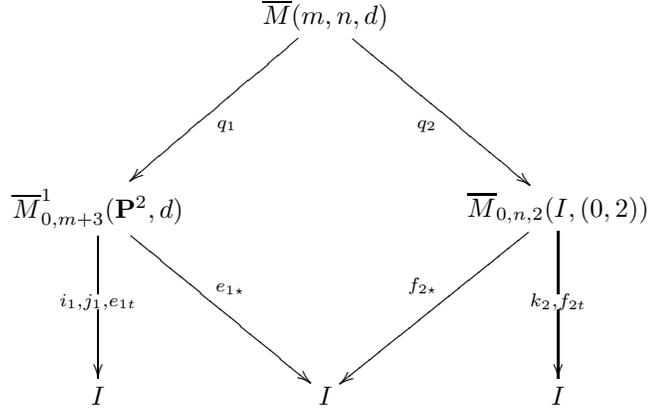
\xymatrix{
& & \spmoduli
\ar[ddl]^{q_1} \ar[ddr]_{q_2} \\
\\
& \lamoduli{m+3}{d} \ar[dd]|{i_1, j_1, e_{1t}} \ar[ddr]^{e_{1\star}} & &
\double{n} \ar[dd]|{k_2, f_{2t}} \ar[ddl]_{f_{2\star}} \\
\\
& I  & I  & I }
$$
Lemma~\ref{pushpull} and the projection formula applied to the above diagram
yield the following expression for the corresponding term of
(\ref{contactproduct}):
\begin{equation}\label{factoredcp}
\frac{2}{m!n!}
k_{2*}\biggl(f_{2\star}^*\bigl(e_{1\star *}
(i_1^*\alpha \cup j_1^*\beta \cup\bigcup_{t=1}^m e_{1t}^*\delta \cap 
[\lamoduli{m+3}{d}])
\bigr)
\cup \bigcup_{t=1}^n f_{2t}^*\delta \cap [\double{n}]
\biggr).
\end{equation}
It follows from Lemma~\ref{lemmapfp} that $\alpha\bullet\beta$ is in $\apt$.
\par
Extending by $\qt$-linearity, we have a product $\bullet$ on $\apt$. We
now define the {\em contact product of $\alpha$ and $\beta$, deformed by $\delta$}, to be 
$$
\alpha * \beta =\alpha\cup\beta +\alpha \bullet \beta.
$$
We call $\apt$, together with this product, the {\em contact cohomology ring} of $\pp$, and denote it by $Q^1H^*(\pp)$.
\begin{thm} For each $\delta$, the contact product is commutative and associative. The identity element $1\in A^0(\pp)$ for the ordinary cup product
also serves as the identity element for the contact product.
\end{thm}
\begin{proof}
The commutativity is obvious. The class
$$
i_1^*1 \cup j_1^*\beta \cup\bigcup_{t=1}^m e_{1t}^*\delta \cap 
[\lamoduli{m+3}{d}]
$$
is the pullback via the forgetful morphism
$$
\lamoduli{m+3}{d}\to \lamoduli{m+2}{d}
$$
of a class on the latter, and $e_{1\star}$ factors through the forgetful
morphism. Hence if $\alpha=1$ then (\ref{factoredcp}) vanishes.
Hence $1$ is the identity for the contact product.
\par
Now, let $\alpha,\beta,\gamma$ be classes in $\ap$. The cup product is associative. Thus we must show that
\begin{equation}
\label{qcpeqn}
(\alpha \cup \beta )\bullet \gamma
+(\alpha \bullet \beta )\cup \gamma
+(\alpha \bullet \beta )\bullet \gamma
=\alpha \cup (\beta \bullet \gamma)
+\alpha \bullet (\beta \cup\gamma)
+\alpha \bullet (\beta \bullet \gamma).
\end{equation}
As in the proof of associativity of the quantum product in section \ref{qtp}, we use a linear equivalence of divisors
$$
D(ij \mid kl) \simeq D(il \mid jk),
$$
this time on the stack of stable lifts $\lamoduli{m+4}{d}$. This linear equivalence, which is derived from the forgetful map $\lamoduli{m+4}{d}\to \stablefour$, induces a linear equivalence on the fiber product $\lamoduli{m+4}{d}\times_I \double{n}$:
$$
D(ij\mid kl)\times_I \double{n}
\simeq D(il \mid jk)\times_I \double{n}.
$$
For these two fiber products, denote the evaluation maps to $I$ coming from the first factor by $i,j,k,l$ and $e_1,\dots,e_m$; denote the evaluation maps coming from the second factor by $r$ and $s$ (for evaluation at the ramification points) and by $f_1,\dots,f_n$. The fiber product is defined via the map $l$ on the first factor and $r$ on the second.
\par
The linear equivalence implies the equality
\begin{multline}
\label{qccdiveqn}
\sum_{\substack{m + n = q\\d>0}}
\frac{2}{m!n!}
s_*\bigl(i^*\alpha \cup j^*\beta \cup k^*\gamma \cup {\bigcup_{t=1}^m}
e_t^*\delta \cup {\bigcup_{t=1}^n}
f_t^*\delta \cap \lbrack D(ij \mid kl)\times_I \double{n} \rbrack \bigr) \\
=
\sum_{\substack{m + n = q\\d>0}}
\frac{2}{m!n!}
s_*\bigl(i^*\alpha \cup j^*\beta \cup k^*\gamma \cup {\bigcup_{t=1}^m}
e_t^*\delta \cup {\bigcup_{t=1}^n}
f_t^*\delta \cap \lbrack D(il \mid jk)\times_I \double{n} \rbrack \bigr).
\end{multline}
We will show that the left side of (\ref{qccdiveqn}) equals the $q$th coefficient of the left side of (\ref{qcpeqn}); an entirely similar argument shows that the right side of (\ref{qccdiveqn}) equals the $q$th coefficient of the right side of (\ref{qcpeqn}).
\par
In \cite[Section 5]{ErnstromKennedy} we have analyzed the components of $D(ij \mid kl)$. We call such a component $D$ {\em numerically irrelevant} if
$$
\int i_1^*\gamma_i\cup j_1^*\gamma_j
\cup k_1^*\gamma_k\cup l_1^*\gamma_l
\cup {\bigcup_{t=1}^m}e_{1t}^*\delta_t
\cap [D]=0.
$$
for all choices of cohomology classes $\gamma_i,\gamma_j,\gamma_k,\gamma_l,$ and $\delta_1, \dots, \delta_m$ in $A^*(I)$. Here $i_1,j_1$, etc. indicate the
maps in the following diagram, for which $i_1 \circ q_1=i$, etc.
\par
$$ \xymatrix{
& & \lamoduli{m+4}{d}\times_I \double{n}
\ar[ddl]^{q_1} \ar[ddr]_{q_2} \\
\\
& \lamoduli{m+4}{d} \ar[dd]|{i_1, j_1, k_1, e_{1t}} \ar[ddr]^{l_1} & &
\double{n} \ar[dd]|{s_2,f_{2t}} \ar[ddl]_{r_2} \\
\\
& I  & I  & I }
$$
\par
To understand the contribution of a numerically irrelevant component to
 (\ref{qccdiveqn}), consider the class
\begin{equation}
\label{vanish}
s_*\bigl(i^*\alpha \cup j^*\beta \cup k^*\gamma \cup {\bigcup_{t=1}^m}
e_t^*\delta \cup {\bigcup_{t=1}^n}
f_t^*\delta \cap \lbrack D \times_I \double{n} \rbrack \bigr).
\end{equation}
Let $\tau$ be a test class in $A_*(I)$.
Then the degree of its intersection with (\ref{vanish})
equals, by repeated use of the projection formula,
the degree of the class
\begin{equation*}
i_1^*\alpha \cup j_1^*\beta \cup k_1^*\gamma \cup {\bigcup_{t=1}^m}
e_{1t}^*\delta 
\cap
q_{1*}\left(
q_2^*(s_2^*\tau \cup{\bigcup_{t=1}^n}
f_{2t}^*\delta)
\cap \lbrack D \times_I \double{n} \rbrack\right).
\end{equation*}
By Lemma~\ref{pushpull}, this class equals
\begin{equation*}
i_1^*\alpha \cup j_1^*\beta \cup k_1^*\gamma \cup {\bigcup_{t=1}^m} e_{1t}^*\delta 
\cup
l_1^*r_{2*}\left(s_2^*\tau \cup {\bigcup_{t=1}^n}f_{2t}^*\delta\cap[\double{n}])\right)
\cap \lbrack D \rbrack.
\end{equation*}
Since $D$ is numerically irrelevant, the degree is 0. And since this is true for every $\tau$, the class (\ref{vanish}) is zero.
\par
Besides the numerically irrelevant components, the divisor $D(ij \mid kl)$ on $\lamoduli{m+4}{d}$ has three types of components \cite[Proposition 5.9]{ErnstromKennedy}. Thus there are three corresponding types of components of $D(ij \mid kl)\times_I \double{n}$. A component of the first type is of the form
\begin{equation}
\label{type1a}
\imoduli{m_1+3}{(0,0)}
\times_I
\lamoduli{m_2+3}{d}
\times_I
\double{n}
\end{equation}
or
\begin{equation}
\label{type1b}
\lamoduli{m_1+3}{d}
\times_I
\imoduli{m_2+3}{(0,0)}
\times_I
\double{n},
\end{equation}
where $m_1+m_2=m$. A typical point of such a component represents a map from a curve with three components, as depicted in Figure 5.2. In each case the component on the right maps two-to-one to a fiber of $I$ and is ramified at the points marked $\times$; the horizontal component maps to $I$ via the lift of an immersion; and the remaining component maps to a point of $I$. In addition to the indicated markings, there are $m_1$ markings on the left component, $m_2$ on the middle component, and $n$ on the right component.
\par
%
%
%
%
$$
\xy
0;<1cm,0cm>:
(0,1);(2.5,1)**\dir{-};
(0.5,1)*{\bullet};
(1,1)*{\bullet}+U(4)*{k};
(2,1)*{\times}*{\bullet}+UL(2.5)*{l}+R(9)*{r};
(0.5,0.75);(0.5,2.5)**\dir{-};
(0.5,1.5)*{\bullet}+L(4)*{j};
(0.5,2)*{\bullet}+L(4)*{i};
(2,0.75);(2,2.5)**\dir{-};
(2,2)*{\times}*{\bullet}+R(4)*{s};
(5,2.25);(7.5,2.25)**\dir{-};
(5.75,2.25)*{\bullet}+U(4)*{i};
(6.5,2.25)*{\bullet}+U(4)*{j};
(7.25,2.5);(7.25,1.5)**\dir{-};
(7.75,1)**\crv{(7.25,1)}*{\times}*{\bullet}+DR(2.5)*{r}+L(9)*{l};
(7.75,1);(8,1)**\dir{-};
(7.75,0.75);(7.75,2.5)**\dir{-};
(7.25,1.875)*{\bullet}+L(4)*{k};
(7.25,2.25)*{\bullet};
(7.75,1.875)*{\bullet}*{\times}+R(2.5)*{s};
\endxy
$$
\begin{center}
{\bf Figure 5.2} Typical members of a component of the first type,
(\ref{type1a}) (left) and (\ref{type1b}) (right.)
\end{center}
\smallskip
\par
Consider a component of type (\ref{type1a}). Let us denote the three factors simply by $M_1,M_2,M_3$. We may relabel the markings so that the evaluation maps $e_1,\dots,e_{m_1}$, as well as $i$, $j$, and the map $g_1$ used to create the fiber product, factor through projection onto $M_1$. We also note that all of these evaluation maps coincide. Hence we have the following fiber diagram.
$$ \xymatrix{
& & M_1\times_I M_2\times_I M_3
\ar[ddl]^{q_{1}} \ar[ddr]_{q_{23}} & \\
\\
& M_1 
\ar[ddr]|{g_{1}=i_{1}=j_{1}=e_{1t}} & &
M_2\times_I M_3  \ar[dd]|{k_{23},e_{23t},f_{23t},s_{23}} 
\ar[ddl]_{g_{23}} & \\
\\
& {} & I  & I  & }
$$
By Lemma~\ref{pushpull} and the projection formula, the contribution of our component to the left side of (\ref{qccdiveqn}) is
\begin{align*}
\frac{2}{m!n!}&
s_*\bigl(i^*\alpha\cup j^*\beta\cup k^*\gamma\cup
\bigcup_{t=1}^me_t^*\delta\cup\bigcup_{t=1}^nf_t^*\delta
\cap [M_1\times_IM_2\times_I M_3]\bigr)\\
&=\frac{2}{m!n!}
s_{23*}\biggl(
g_{23}^*g_{1*}\biggl(
g_{1}^*(\alpha\cup \beta \cup \bigcup_{t=1}^{m_1}\delta) 
\cap [M_1] \biggr)
\cup
k_{23}^*\gamma
\cup
\bigcup_{t=m_1+1}^{m}e_{23t}^*\delta
\cup
\bigcup_{t=1}^n f_{23t}^*\delta
\cap [M_2\times_IM_3]
\biggr).
\end{align*}
If $m_1 > 0$ then the fibers of the map from $M_1$ to $I$ have positive dimension; hence  
\begin{equation*}
g_{1*}\biggl(
g_{1}^*(\alpha\cup \beta \cup \bigcup_{t=1}^{m_1}\delta) 
\cap [M_1] \biggr) = 
\alpha\cup\beta\cup\bigcup_{t=1}^{m_1}\delta\cap g_{1*}[M_1]=0,
\end{equation*}
and the component makes no contribution to (\ref{qccdiveqn}). If $m_1 = 0$ then $M_1$ and $I$ are isomorphic, and
\begin{equation*}
g_{23}^*g_{1*}\biggl(
g_{1}^*(\alpha\cup \beta \cup \bigcup_{t=1}^{m_1}\delta) 
\cap [M_1] \biggr) = 
g_{23}^*\biggl(\alpha\cup \beta 
\cap [I] \biggr).
\end{equation*}
Thus if we sum the contributions from all such components with $m+n=q$ we obtain the $q$th coefficient of $(\alpha \cup \beta )\bullet \gamma$.
\par
Entirely similar arguments show that a component of type (\ref{type1b}) makes no contribution to (\ref{qccdiveqn}) unless $m_2=0$, and that the sum of contributions from all such components with $m+n=q$ is the $q$th coefficient of $(\alpha \bullet \beta )\cup \gamma$.
\par
The second type of component of $D(ij \mid kl)\times_I \double{n}$ is one of the form
\begin{equation}\label{type2a}
{\overline M}_{0,m_1+3}(I,(0,1))
\times_I {\overline C}^1_{0,m_2+2,\{\star\}}(\pp,d)
\times_I \double{n},
\end{equation}
or
\begin{equation}\label{type2b}
{\overline C}^1_{0,m_1+2,\{\star\}}(\pp,d)
\times_I
{\overline M}_{0,m_2+3}(I,(0,1))
\times_I \double{n}
\end{equation}
where $m_1+m_2=m$, and where the second (respectively first) factor is the stack of cuspidal stable lifts \cite[Section 4]{ErnstromKennedy}. A general point of such a component 
represents a map from a curve with the same configuration of components and
markings as in Figure 5.2. The left and middle components map,
in either order, to the lift of a degree $d$
rational curve with one cusp, and to the
fiber of $I$ over the cusp, and the right component maps two-to-one to a
fiber of $I$.
\par
Consider a component of type (\ref{type2a}). Denote the three
factors by $M_1, M_2, M_3$ and consider the following fiber diagram.
\par
$$ \xymatrix{
 & M_1\times_I M_2 \times_I M_3
\ar[ddl]^{q_1} \ar[ddr]_{q_{23}} \\
\\
 M_1 \ar[dd]|{i_1, j_1, e_{1t}} \ar[ddr]^{g_1} & &
 M_2 \times_I M_3 \ar[dd]|{k_{23},s_{23},e_{23t},f_{23t}} \ar[ddl]_{g_{23}} \\
\\
 I  & I  & I }
$$
By Lemma~\ref{pushpull} and the projection formula, the contribution of our
component to the left side of (\ref{qccdiveqn}) is
$$
\frac{2}{m! n!}s_{23*}\biggl(g_{23}^*g_{1*}\bigl(i_1^*\alpha\cup
j_1^*\beta\cup\bigcup_{t=1}^{m_1}e_{1t}^*\delta \cap [M_1]\bigr)
\cup k_{23}^*\gamma\cup\bigcup_{t=m_1+1}^m e_{23t}^* \delta\cup
\bigcup_{t=1}^nf_{23t}^*\delta\cap [M_2\times_I M_3]\biggr).
$$
By Lemma~\ref{disappear}, the class
$$
g_{1*}\bigl(i_1^*\alpha\cup
j_1^*\beta\cup\bigcup_{t=1}^{m_1}e_{1t}^*\delta \cap [M_1]\bigr)=0.
$$
Hence the component makes no contribution.
A similar argument applies to a component of type (\ref{type2b}).
\par
We come finally to the third type of component of $D(ij \mid kl)\times_I \double{n}$. Among them are components of the form
\begin{equation}
\label{type3a}
\lamoduli{m_1+3}{d_1}\times_I \double{m_2}
\times_I\lamoduli{m_3+3}{d_3}\times_I\double{n},
\end{equation}
with $m_1+m_2+m_3=m$, $d_1,d_3>0$ and $d_1+d_3=d$.
A typical point of such a component represents a map from a curve with four components, as depicted at the top of Figure 5.3. The vertical components map two-to-one to fibers of $I$ and are ramified at the points marked $\times$; the horizontal components map to $I$ via the lifts of immersions. In addition to the indicated markings, there are $m+n$ others.
\par
%
%
%
%
$$
\xy
0;<1cm,0cm>:
(0,1.75);(1.5,1.75)**\dir{-};
(0.25,1.75)*{\bullet}+U(4)*{i};
(0.75,1.75)*{\bullet}+U(4)*{j};
(1.25,2);(1.25,0)**\dir{-};
(1.25,0.25)*{\times}*{\bullet};
(1.25,1.75)*{\times}*{\bullet};
(1.25,0.75)*{\bullet}+L(4)*{k};
(1,0.25);(3,0.25)**\dir{-};
(2.75,0.25)*{\times}*{\bullet}+UL(2.5)*{l}+R(9)*{r};
(2.75,0);(2.75,1.75)**\dir{-};
(2.75,1.5)*{\times}*{\bullet}+R(4)*{s};
(0,4.75);(1.5,4.75)**\dir{-};
(0.25,4.75)*{\bullet}+U(4)*{i};
(1.25,5);(1.25,3)**\dir{-};
(1.25,3.25)*{\times}*{\bullet};
(1.25,4.75)*{\times}*{\bullet};
(1.25,3.75)*{\bullet}+L(4)*{j};
(1,3.25);(3,3.25)**\dir{-};
(2,3.25)*{\bullet}+U(4)*{k};
(2.75,3.25)*{\times}*{\bullet}+UL(2.5)*{l}+R(9)*{r};
(2.75,3);(2.75,4.75)**\dir{-};
(2.75,4.5)*{\times}*{\bullet}+R(4)*{s};
(4,4.75);(5.5,4.75)**\dir{-};
(5.25,5);(5.25,3)**\dir{-};
(5.25,3.25)*{\times}*{\bullet};
(5.25,4.75)*{\times}*{\bullet};
(5.25,3.75)*{\bullet}+L(4)*{j};
(5.25,4.25)*{\bullet}+L(4)*{i};
(5,3.25);(7,3.25)**\dir{-};
(6,3.25)*{\bullet}+U(4)*{k};
(6.75,3.25)*{\times}*{\bullet}+UL(2.5)*{l}+R(9)*{r};
(6.75,3);(6.75,4.75)**\dir{-};
(6.75,4.5)*{\times}*{\bullet}+R(4)*{s};
(8,4.75);(9.5,4.75)**\dir{-};
(8.25,4.75)*{\bullet}+U(4)*{j};
(9.25,5);(9.25,3)**\dir{-};
(9.25,3.25)*{\times}*{\bullet};
(9.25,4.75)*{\times}*{\bullet};
(9.25,3.75)*{\bullet}+L(4)*{i};
(9,3.25);(11,3.25)**\dir{-};
(10,3.25)*{\bullet}+U(4)*{k};
(10.75,3.25)*{\times}*{\bullet}+UL(2.5)*{l}+R(9)*{r};
(10.75,3);(10.75,4.75)**\dir{-};
(10.75,4.5)*{\times}*{\bullet}+R(4)*{s};
(4,7.75);(5.5,7.75)**\dir{-};
(4.25,7.75)*{\bullet}+U(4)*{i};
(4.75,7.75)*{\bullet}+U(4)*{j};
(5.25,8);(5.25,6)**\dir{-};
(5.25,6.25)*{\times}*{\bullet};
(5.25,7.75)*{\times}*{\bullet};
(6,6.25)*{\bullet}+U(4)*{k};
(5,6.25);(7,6.25)**\dir{-};
(6.75,6.25)*{\times}*{\bullet}+UL(2.5)*{l}+R(9)*{r};
(6.75,6);(6.75,7.75)**\dir{-};
(6.75,7.5)*{\times}*{\bullet}+R(4)*{s};
(4,1.75);(6,1.75)**\dir{-};
(4.25,1.75)*{\bullet}+U(4)*{i};
(4.75,1.75)*{\bullet}+U(4)*{j};
(5.75,1.75)*{\times}*{\bullet};
(5.75,2);(5.75,0)**\dir{-};
(5.75,1)*{\times}*{\bullet};
(5.5,1);(7,1)**\dir{-};
(5.25,1.5);(5.25,1)**\dir{-}*{\times}*{\bullet}+L(4)*{s};
(5.25,1);(5.75,0.5)**\crv{(5.25,0.5)};
(6,0.5);(5.75,0.5)**\dir{-}*{\times}*{\bullet}+DR(2.5)*{l}+L(9)*{r};
(5.75,1.4)*{\bullet}+R(4)*{k};
(8,1.75);(10,1.75)**\dir{-};
(8.25,1.75)*{\bullet}+U(4)*{i};
(8.75,1.75)*{\bullet}+U(4)*{j};
(9.75,1.75)*{\times}*{\bullet};
(9.75,2);(9.75,0)**\dir{-};
(9.75,1)*{\times}*{\bullet};
(9.5,1);(11,1)**\dir{-};
(9.75,1.75)*{\bullet};
(9.25,1.5);(9.25,1)**\dir{-}*{\times}*{\bullet}+L(4)*{s};
(9.25,1);(9.75,0.5)**\crv{(9.25,0.5)};
(10,0.5);(9.75,0.5)**\dir{-}*{\times}*{\bullet}+DR(2.5)*{l}+L(9)*{r};
(10.3625,1)*{\bullet}+U(4)*{k}
\endxy
$$
\begin{center}
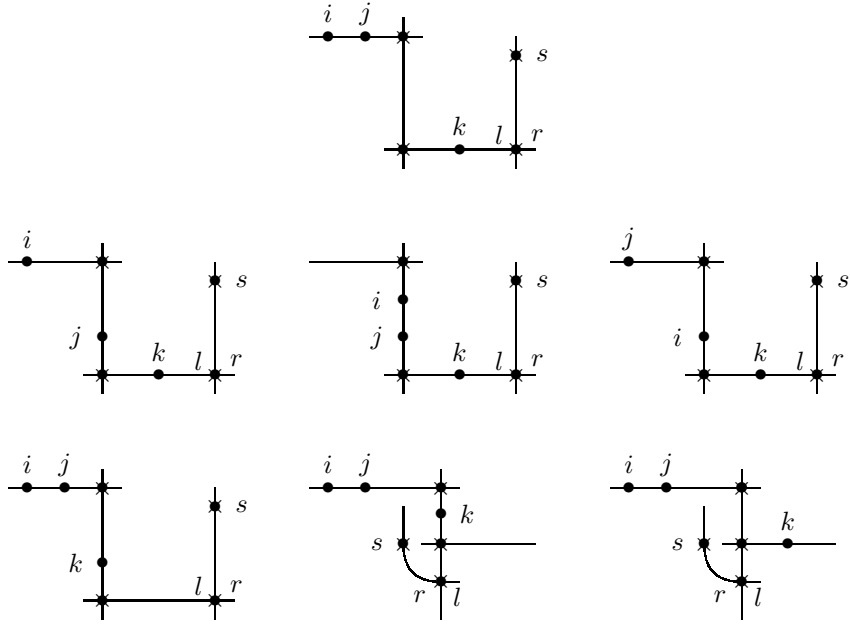
{\bf Figure 5.3} The different configurations for components of the third type.
\end{center}
\smallskip
\par
There are six other possibilities for a component of the third type, corresponding to the six different possible ways in which the four special markings $i,j,k,l$ can lie on the first three components of the typical curve. (The markings $i$ and $j$ may lie on either of the first two components; $k$ and $l$ may lie on either the second or third component; and the pair $i,j$ must be separated from the pair $k,l$ by a node.) We claim that in each of these six cases the component makes no contribution to (\ref{qccdiveqn}). In five cases the argument is identical to that for a component of the second type: a configuration with $i,j$ or $k$ on a component mapping to a fiber (vertical component in Figure 5.3) does not contribute, as a consequence of Lemma~\ref{lemmapbpp}.

In the sixth case we have a component of the form
\begin{equation}
\label{type3f}
\biggl(\lamoduli{m_1+3}{d_1}\times_I \double{m_2+1}
\times_I\lamoduli{m_3+2}{d_3}\biggr)\times_I\double{n},
\end{equation}
A typical point of such a component represents a map from the sort of
curve shown at the bottom right of Figure 5.3. Denote the four factors
by $M_1, M_2, M_3$ and $M_4$. Note that the fiber product of
$M_1\times_I M_2 \times_I M_3$ and $M_4$ is created by using an evaluation
map from $M_1\times_I M_2 \times_I M_3$ which comes from its second factor.
Consider the following fiber diagram.
\par
$$ 
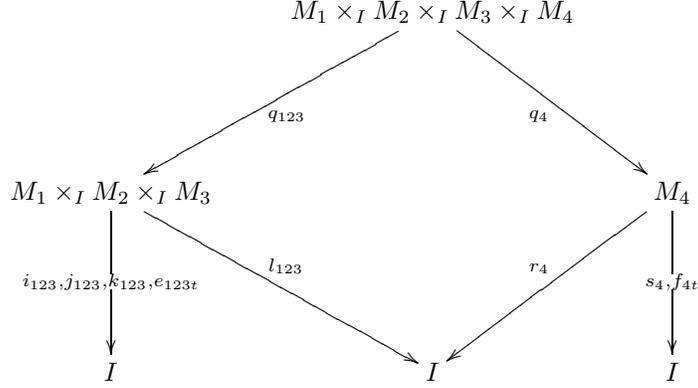
\xymatrix{
 & M_1\times_I M_2 \times_I M_3\times_I M_4
\ar[ddl]^{q_{123}} \ar[ddr]_{q_4} \\
\\
 M_1\times_I M_2 \times_I M_3  \ar[dd]|{i_{123}, j_{123}, k_{123}, e_{123t}} \ar[ddr]^{l_{123}} & &
M_4 \ar[dd]|{s_4,f_{4t}} \ar[ddl]_{r_4} \\
\\
 I  & I  & I }
$$
Using Lemma~\ref{pushpull} one finds that the contribution to (\ref{qccdiveqn})
is
$$
\frac{2}{m!n!}s_{4*}\biggl(\bigcup_{t=1}^{n}f_{4t}^*\delta
\cup
r_4^* l_{123*}\bigl(i_{123}^*\alpha\cup j_{123}^*\beta\cup k_{123}^*\gamma \cup
\bigcup_{t=1}^{m}e_{123t}^*\delta
\cap [M_1\times_I M_2 \times_I M_3]\bigr)
\cap [M_4]\biggr).
$$
We claim that the class
\begin{equation}\label{ppclass}
l_{{123}*}\bigl(i_{123}^*\alpha\cup j_{123}^*\beta\cup k_{123}^*\gamma
\cup\bigcup_{t=1}^{m}e_{{123}t}^*\delta
\cap [M_1\times_I M_2 \times_I M_3]\bigr)
\end{equation}
is in $\ap$. It will then follow from Lemma~\ref{lemmapbpp} that the contribution is zero.
To prove the claim we use the following diagram, in which there are three fiber squares.
$$ 
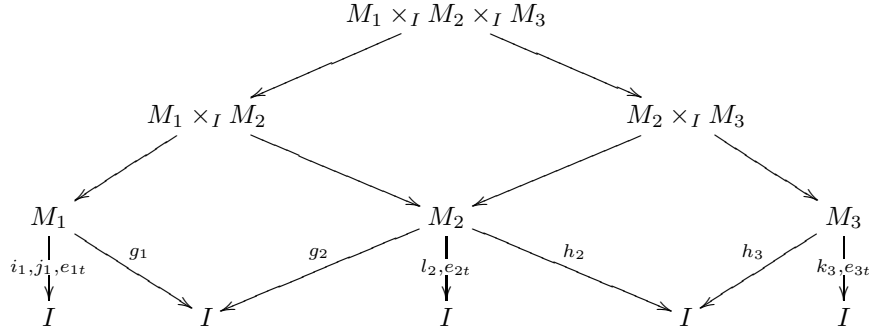
\xymatrix{
 & &  M_1\times_I M_2 \times_I M_3
\ar[dl]  \ar[dr] \\
 & M_1\times_I M_2 \ar[dl]  \ar[dr] 
& & M_2\times_I M_3 \ar[dl] \ar[dr]\\
 M_1 \ar[d]|{i_{1}, j_{1}, e_{{1}t}} \ar[dr]^{g_1} & & 
M_2 \ar[dl]_{g_2} \ar[d]|{l_2, e_{2t}} \ar[dr]^{h_2} & &
M_3 \ar[dl]_{h_3} \ar[d]|{k_3,e_{3t}}\\
I & I & I & I & I }
$$
By repeated use of Lemma~\ref{pushpull} and the projection formula, the class
(\ref{ppclass}) can be expressed as
$$
l_{2*}\biggl(
g_2^*g_{1*}\bigl(i_1^*\alpha\cup j_1^*\beta\cup
\bigcup_{t=1}^{m_1}e_{1t}^*\delta\cap [M_1]\bigr)
\cup h_2^* h_{3*}\bigl(k_3^*\gamma\cup
\bigcup_{t=m_1+m_2+1}^{m_1+m_2+m_3}e_{3t}^*\delta\cap [M_3]\bigr)
\cup\bigcup_{t=m_1+1}^{m_1+m_2}e_{2t}^*\delta\cap [M_2]\biggr).
$$
Now $M_2=\double{n}$, so by Lemma~\ref{lemmapfp} this is a class
in $\ap$ as we claimed.
\par
It remains to consider the contribution of the component with the
configuration on the top of Figure 5.3. Note that the product (\ref{type3a})
can be written as the following fiber product:
\par
$$ 
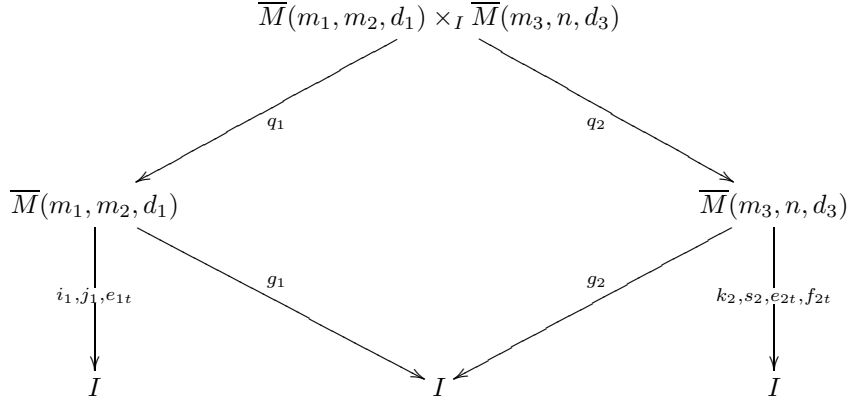
\xymatrix{
& & \spamoduli{m_1}{m_2}{d_1}\times_I\spamoduli{m_3}{n}{d_3}
\ar[ddl]^{q_1} \ar[ddr]_{q_2} \\
\\
&  \spamoduli{m_1}{m_2}{d_1} \ar[dd]|{i_1, j_1, e_{1t}} \ar[ddr]^{g_1} & &
\spamoduli{m_3}{n}{d_3} \ar[dd]|{k_2,s_2,e_{2t},f_{2t}} \ar[ddl]_{g_2} \\
\\
& I  & I  & I }
$$
Using Lemma~\ref{pushpull} one finds that the contribution to the left side of (\ref{qccdiveqn}) is
\begin{multline*}
\frac{2}{n!m!}s_{2*}\biggl(k_2^*\gamma\cup\bigcup_{t=m_1+m_2+1}^{m_3}e_{2t}^*\delta
\cup\bigcup_{t=1}^n f_{2t}^*\delta\cup
g_2^*g_{1*}\bigl(i_1^*\alpha\cup j_1^*\beta\cup \\
\bigcup_{t=1}^{m_1}e_{1t}^*\delta\cup\bigcup_{t=m_1+1}^{m_2}e_{1t}^*\delta
\cap [\spamoduli{m_1}{m_2}{d_1}]\bigr)
\cap [\spamoduli{m_3}{n}{d_3}]\biggr).
\end{multline*}
According to \cite{ErnstromKennedy} each divisor component of this type
appears in $D(ij \mid kl)$ with a multiplicity of two. The reason is that 
$(d,2d-2)$ can be partitioned in two ways, either
as $(d_1,2d_1-2)+(d_2,2d_2)$ or as
$(d_1,2d_1)+(d_2,2d_2-2)$.
If we keep  $q$ fixed,
and sum the contributions for all $d$,
$d_1+d_3=d$, $m+n=q$ and all partitions of the $m$ markings into three
disjoint subsets with cardinalities $m_1$, $m_2$ and $m_3$,
the result is the $q$th coefficient of
$(\alpha\bullet\beta)\bullet\gamma$.

Thus we have shown that $(\alpha * \beta)*\gamma=\alpha *(\beta * \gamma)$
when $\alpha, \beta, \gamma$ are elements of $\ap$. Invoking the
$\qt$-linearity, we conclude that the contact product is associative.
\par
\end{proof}
\section{The recursive relation among characteristic numbers}
\label{rrc}
The associativity of the quantum product implies Kontsevich's recursive formula for the characteristic numbers
\begin{equation*}
N_d=
\text{ the number of rational plane curves of degree $d$
through $3d-1$ general points.}
\end{equation*}
For details of this story, see \cite{FrancescoItzykson}, \cite{FultonP}, \cite{KontsevichManin}. Here we show that, in a similar fashion, the associativity of the contact product implies a recursive formula for the numbers 
\begin{align*}
N_d(a,b,c)=&\text{ the number of rational plane curves of degree $d$
through $a$ general points, tangent to}\\ &\text{$b$ general lines, and
tangent to $c$ general lines at a specified general point on each line}\\
&\quad\text{(where $a+b+2c=3d-1$).}
\end{align*}
Our formula will specialize both to Kontsevich's formula and the more general formula of Di~Francesco and Itzykson \cite[Equation 2.95]{FrancescoItzykson}.
\par
We will use the following ordered basis for the cohomology of the incidence correspondence:
$$
\{T_0,T_1,T_2,T_3,T_4,T_5\}=\{1,h,h^2,\hd,\hd^2,h^2\hd\}.
$$
With respect to this basis the fundamental class of the diagonal $\Delta$
in $I \times I$ has the simple decomposition
\begin{equation*}
[\Delta] = \sum_{s=0}^5 [T_s] \times [T_{5-s}].
\end{equation*}
Suppose that $d$ is a positive integer, and that
$\gamma_1,\dots,\gamma_n$ are elements of $A^*I$.
Then the {\em first-order Gromov-Witten invariant} is
$$
N_d(\gamma_1\cdots\gamma_n)=\int e_1^*(\gamma_1) \cup \dots \cup
e_n^*(\gamma_n)
\cap [\lamoduli{n}{d}].
$$
According to \cite[Section 4]{ErnstromKennedy}, we have the following interpretations:
\begin{trivlist}
\item[(1)] Suppose that $a$ of the $\gamma_t$'s equal the class $h^2$, that $b$ of them equal the class $\hd^2$, and that the remaining $c$ of them equal $h^2\hd$, where $a+b+2c=3d-1$. Then the Gromov-Witten invariant is the characteristic number $N_d(a,b,c)$.
\item[(2)] For all $d$ and all $\gamma_1,\dots,\gamma_{n-1}$,
$$
N_d(\gamma_1\cdots\gamma_{n-1}\cdot 1)=0. $$
\par
\item[(3)] If $\gamma_n$ is the class of a divisor, then for all $d$ and all $\gamma_1,\dots,\gamma_{n-1}$,
$$
N_d(\gamma_1\cdots\gamma_n)=
N_d(\gamma_1\cdots\gamma_{n-1})\int\gamma_n\cap [C], $$
where $[C]=d \hd^2 +(2d-2) h^2$.
\end{trivlist}
\par
Using a general element
$$
\gamma=y_0T_0+\dots+y_5T_5
$$
of $A^*I$, we define the {\em quantum potential} to be the following formal power series in $y_0,\dots,y_5$:
$$
\calN=\sum_{\substack{m \geq 0 \\ d\geq 1}}  \frac{1}{m!}
\int e^*_1(\gamma)\cup \dots \cup e^*_m(\gamma)
\cap [\lamoduli{m}{d}],
$$
where $e_1,\dots,e_m$ are the evaluation maps. By the previous remarks
$$
\calN=\sum_{\substack{d\geq 1 \\a +b+2c=3d-1\\ a,b,c\geq 0}}
\frac{N_d(a,b,c)y_2^ay_4^by_5^c\exp(dy_1+(2d-2)y_3)}{a!b!c!}.
$$
In a similar way, we define a potential associated to the stacks
$\double{n}$. This will be a formal power series in two sets of indeterminates. Let $\delta=z_0T_0+\dots+z_5T_5$ be a second general element of $A^*I$. Then
\begin{equation*}
\calR=\sum_{n \geq 0} \frac{1}{2n!}
\int e^*_1(\gamma) \cup \dots \cup e^*_n(\gamma)
\cup e^*_{n+1}(\delta) \cup  e^*_{n+2}(\delta)
\cap [\double{n}],
\end{equation*}
where $e_{n+1}$ and $e_{n+2}$ are evaluation at the points of ramification. As we show in \cite[Section 6]{ErnstromKennedy},
\begin{equation}
\label{potr}
\calR=\left\{\frac{z_3^2}{2}(y_4^2+y_5)+
z_3z_4y_4+\frac{z_3z_5}{2}+\frac{z_4^2}{4}\right\}\exp(2y_3).
\end{equation}
\par
Similarly, we define a potential associated to the stacks
$$
\lamoduli{m+1}{d} \times_I \double{n}.
$$
Let $e_1,\dots,e_m$ and $e_{\star}$ be the evaluation maps coming from the first factor; let $f_{\star}$ and $k$ be evaluation at the points of ramification; let $f_1,\dots,f_n$ be the other evaluations coming from the second factor; let $e_{\star}$ and $f_{\star}$ be the maps used to create the fiber product. We define $\calK$ to be
\begin{equation*}
\sum_{\substack{m, n \geq 0 \\ d\geq 1}}
\frac{2}{m!n!}
\int e^*_1(\gamma) \cup \dots \cup e^*_m(\gamma)
\cup f^*_1(\gamma) \cup  \dots \cup f^*_n(\gamma)
\cup k^*(\delta)
\cap
[\lamoduli{m+1}{d} \times_I \double{n}].
\end{equation*}
Then by Propositions 6.2 and 6.3 of \cite{ErnstromKennedy}, the three potentials are related by the differential equation
\begin{equation}
\label{knr}
\calK= 2 \sum_{s=0}^5
\frac{\partial \calN}{\partial y_s}
\frac{\partial \calR}{\partial z_{5-s}};
\end{equation}
from the explicit form (\ref{potr}) of $\calR$, however, we see that the sum needs to run only from $s=0$ to $2$.
\par
Proposition 6.3 also tells us how the product $\bullet$ is related to the potential $\calK$: for $0 \leq i,j \leq 2$ we have
$$
T_i \bullet T_j = \sum_{s=0}^5
\kp{i}{j}{5-s} T_s.
$$
Again we note that the sum needs to run only from $s=0$ to $2$, since the product of two elements of $\ap$ is a formal power series whose coefficients are likewise in $\ap$. Thus
$$
T_i \bullet T_j = 
\kp{i}{j}{5} T_0 +
\kp{i}{j}{4} T_1 +
\kp{i}{j}{3} T_2.
$$
\par
We have shown that the contact product is associative. In particular $(T_1*T_1) * T_2 = T_1*(T_1*T_2)$. Equating the coefficients of $T_0$ on the two sides of this equation, we find that
\begin{multline*}
\kp{1}{1}{4}\kp{1}{2}{5}+\kp{2}{2}{5}+\kp{1}{1}{3}\kp{2}{2}{5} \\
=\kp{1}{1}{5}\kp{1}{2}{4}+\kp{1}{2}{5}\kp{1}{2}{3}.
\end{multline*}
Applying (\ref{knr}) throughout and simplifying, we obtain the following
partial differential equation for the potential $\calN$:
\begin{equation}
\label{npde}
\calN_{222}=\exp(2y_3)\biggl(\calN_{112}^2-\calN_{111}\calN_{122}
+2y_4\bigl(\calN_{112}\calN_{122}-\calN_{111}\calN_{222}\bigr)
+(2y_4^2+2y_5)\bigl(\calN_{122}^2-\calN_{112}\calN_{222}\bigr)\biggr).
\end{equation}
Note that if we set $y_3=y_5=0$ we recover equation (2.95) of \cite{FrancescoItzykson}, and that if furthermore we set $y_4=0$ we recover equation (5.16) of \cite{KontsevichManin}.
\par
Our equation (\ref{npde}) can be rewritten as a formula for the characteristic
number $N_d(a,b,c)$. Note we must assume that $a\geq 3$. In this formula $d_1$ and $d_2$ are greater than zero; thus it determines our characteristic number if we assume we already know those for curves of lower degree.
\begin{equation}
\label{recurs}
\begin{aligned}
N_d(a,b&,c)=\sum_{\substack{d_1+d_2=d\\ a_1+a_2=a-1\\b_1+b_2=b\\c_1+c_2=c}}
N_{d_1}(a_1,b_1,c_1)N_{d_2}(a_2,b_2,c_2)
\left[d_1^2d_2^2\binom{a-3}{a_1-1}
-d_1^3d_2\binom{a-3}{a_1}
\right]\binom b{b_1}\binom c{c_1}\\
&+2\cdot\sum_{\substack{d_1+d_2=d\\ a_1+a_2=a\\b_1+b_2=b-1\\c_1+c_2=c}}
N_{d_1}(a_1,b_1,c_1)N_{d_2}(a_2,b_2,c_2)
\left[d_1^2d_2\binom{a-3}{a_1-1}
-d_1^3\binom{a-3}{a_1}
\right]\binom b{b_1\, b_2\, 1}\binom c{c_1}\\
&+4\cdot\sum_{\substack{d_1+d_2=d\\ a_1+a_2=a+1\\b_1+b_2=b-2\\c_1+c_2=c}}
N_{d_1}(a_1,b_1,c_1)N_{d_2}(a_2,b_2,c_2)
\left[d_1d_2\binom{a-3}{a_1-2}
-d_1^2\binom{a-3}{a_1-1}
\right]\binom b{b_1\, b_2\, 2}\binom c{c_1}\\
&+2\cdot\sum_{\substack{ d_1+d_2=d\\ a_1+a_2=a+1\\b_1+b_2=b\\c_1+c_2=c-1}}
N_{d_1}(a_1,b_1,c_1)N_{d_2}(a_2,b_2,c_2)
\left[d_1d_2\binom{a-3}{a_1-2}
-d_1^2\binom{a-3}{a_1-1}
\right]\binom b{b_1}\binom c{c_1\, c_2\, 1}.
\end{aligned}
\end{equation}
If we set $c=0$ we recover the recursive formula of
Di~Francesco and Itzykson \cite{FrancescoItzykson}. They never state the formula explicitly, but it is an immediate consequence of their equation (2.95), and they have clearly used it in calculating the table of values (2.97).
If we also set $a=3d-1$ and $b=0$ we recover the formula (5.17) of \cite[Claim 5.2.1]{KontsevichManin}.
\par
Formula \ref{recurs} by itself is not enough to determine all characteristic numbers, since we need to supply it with all those values for which $a \leq 2$. Now it may appear that we could perhaps derive additional information by writing down the associativity equations for three other basis elements (rather than $T_1, T_1, T_2$) or by extracting the coefficients with respect to some other basis element (rather than $T_0$). But in fact all such derivations lead---if not to a triviality---to the same differential equation (\ref{npde}). We will omit the proof of this fact. To determine all characteristic numbers it is necessary to use the other relations of \cite[Section 7]{ErnstromKennedy}.
\par
Using the formalism of Fulton and Pandharipande \cite[9 Prop.10]{FultonP}, it is possible to give a presentation of the contact cohomology ring, based on the standard presentation
$$
0\to (z^3)\to \q[z]\to A^*(\pp)\to 0
$$
of $\ap$ as an algebra over $\q$, in which $z^i$ is sent to the class $h^i=h\cup\dots\cup h$. 
Using their arguments, we see that $1,h$ and $h*h$ form a
$\q[[y_0,\dots,y_5]]$-basis for $Q^1H^*(\pp)$, and that we have a similar presentation for the $\q[[y_0,\dots,y_5]]$-algebra $Q^1H^*(\pp)$:
$$
0\to (z^3-\xi_2z^2-\xi_1z-\xi_0)\to\q[[y_0,\dots,y_5]][z]\to Q^1H^*(\pp)\to 0.
$$
Here $z^i$ is sent to $h^{*i}=h * \dots * h$, and thus the elements 
$\xi_i$ are the coefficients of $h*h*h$ with respect to this basis.
Explicitly,
\begin{multline}
Q^1H^*(\pp)=\q[[y_0,\dots,y_5]][z]/\biggl(z^3-\exp(2y_3)\biggl(
\bigl(\calN_{111}+4y_4\calN_{112}+(2y_4^2+2y_5)\calN_{122}\bigr)z^2 \\
+\bigl(2\calN_{112}+2y_4\calN_{122}\bigr)z+\calN_{122}\biggr)
-\exp(4y_3)\bigl(\calN_{111}\calN_{122}-\calN_{112}^2\bigr)
\bigl((2y_4^2+2y_5)z+2y_4\bigr)\biggr).
\end{multline}
If we set $y_3=y_4=y_5=0$ we recover their presentation of $QH^*(\pp)$
\cite[9 Eqn.64]{FultonP}.
\bibliographystyle{nyjalpha}
\ifx\undefined\bysame
\newcommand{\bysame}{\leavevmode\hbox to3em{\hrulefill}\,}
\fi

\end{document}